\documentclass[11pt,reqno]{amsart}

\usepackage{amsmath, amsfonts, amssymb, amscd, enumerate}
\theoremstyle{plain}
\newtheorem{theo}{Theorem}[section]
\newtheorem{lem}[theo]{Lemma}
\newtheorem{prop}[theo]{Proposition}

\theoremstyle{definition}
\newtheorem{rem}[theo]{Remark}

\newtheorem{definition}[theo]{Definition}
\newenvironment{coord}{\noindent{\it Expressions in coordinates. }}{\par\medskip}
\newenvironment{pf}{\noindent{\it Proof. }}{$\square$\par\medskip}

\theoremstyle{plain}

\theoremstyle{definition}

\renewcommand{\=}{\overset{\operatorname{def}}{=}}

\newcommand{\beq}{\begin{equation}}
\newcommand{\eeq}{\end{equation}}
\renewcommand{\a}{\alpha}

\renewcommand{\d}{\delta}

\newcommand{\g}{\gamma}

\renewcommand{\l}{\lambda}
\renewcommand{\o}{\omega}

\renewcommand{\r}{\rho}
\newcommand{\s}{\sigma}

\renewcommand{\O}{\Omega}

\newcommand{\bR}{\mathbb{R}}

\newcommand{\gd}{\mathfrak{d}}

\newcommand{\gF}{\mathfrak{F}}
\newcommand{\gX}{\mathfrak{X}}

\newcommand{\so}{\mathfrak{so}}

\newcommand\GL{\mathrm{GL}}

\newcommand\Spin{\mathrm{Spin}}

\newcommand{\cG}{\mathcal{G}}
\newcommand{\cH}{\mathcal{H}}

\newcommand{\cL}{\mathcal{L}}

\newcommand{\cP}{\mathcal{P}}

\newcommand{\cS}{\mathcal{S}}
\newcommand{\cT}{\mathcal{T}}
\newcommand{\cU}{\mathcal{U}}

\newcommand\Gam[3]{\Gamma^{\ #1}_{\!\!#2#3}}

\renewcommand{\square}{\kern1pt\vbox
{\hrule height 0.6pt\hbox{\vrule width 0.6pt\hskip 3pt
\vbox{\vskip 6pt}\hskip 3pt\vrule width 0.6pt}\hrule height0.6pt}\kern1pt}

\newcommand{\Aff}{{\operatorname{Aff}}}
\newcommand{\Gau}{\mathrm{Gau}}

\newcommand{\Conn}{{\mathrm{Conn}}}
\newcommand{\Grav}{{\mathrm{Grav}}}

\newcommand{\reg}{{\mathrm{reg}}}

\newcommand{\wt}{\widetilde}
\newcommand{\wh}{\widehat}

\newcommand{\n}{\nabla}

\newcommand{\be}{\begin{equation}}
\newcommand{\ee}{\end{equation}}

\def\<#1,#2>{\langle\,#1,\,#2\,\rangle}
\newcommand{\arr}{\begin{array}{rlll}}
\newcommand{\ea}{\end{array}}
\newcommand{\bea}{\begin{eqnarray}}
\newcommand{\eea}{\end{eqnarray}}
\newcommand{\bean}{\begin{eqnarray*}}
\newcommand{\eean}{\end{eqnarray*}}

\def\sideremark#1{\ifvmode\leavevmode\fi\vadjust{
\vbox to0pt{\hbox to 0pt{\hskip\hsize\hskip1em
\vbox{\hsize3cm\tiny\raggedright\pretolerance10000
\noindent #1\hfill}\hss}\vbox to8pt{\vfil}\vss}}}

\newcounter{ssig}
\newcounter{ttig}
\title[The two ways of  gauging the Poincar\`e group]
{The two ways of\\
gauging the Poincar\`e group}

\author{A. Spiro and S. Tantucci}

\begin{document}

\begin{abstract} A   description of how a theory of gravity can be considered as a gauge theory  (in the sense of  Trautman) of  the Poincar\`e group is given. As a result, it is shown that  a gauge theory of this kind is  consistent  with  the Equivalence Principle only if  the Lagrangian and the constraints are  preserved not only by  the gauge transformations  but also by an additional  family of transformations, called {\it pseudo-translations}. Explicit  expressions of pseudo-translations and of their action on gravitational gauge fields are given. They   are expected to be useful for 
 geometric interpretations of their analogues   in  supergravity theories. 
\end{abstract}

\address{
Dipartimento di Matematica e Informatica, Universit\`a di Camerino, 
 Camerino (MC), 
ITALY}
\email{andrea.spiro@unicam.it}

\address{
Via Goggetta n. 13, 
60024 Filottrano (AN),
ITALY}
\email{silvano.tantucci@gmail.com}
\maketitle
\section{Introduction}
\setcounter{section}{1}
\setcounter{equation}{0}

The notion of  ``gauge theory''  is certainly very well-known, but we need to explicitly  recall it in order to make clear  in which sense we are going to use it.  Following Trautman (\cite{Tr, Tr1}), we call {\it gauge theory of a Lie group $G$\/}  a triple formed by:
\begin{itemize}
\item[a)] a class of principal $G$-bundles $\pi: P\to M$ over a manifold $M$ and endowed with connections; 
the manifold $M$ represents the physical empty space-time, 
while the connections on $P$  represent (potentials for) physical fields, called ``gauge fields''; 
\item[b)] a class of bundles $\pi^E: E \to M$,  associated to the previous  $G$-bundles, whose  sections  represent  other physical   fields,   called  ``particles coupled with the gauge fields''; 
\item[c)] a set of partial differential equations (usually Euler-Lagrange equations determined by a Lagrangian) on the connections in (a) and  the sections in (b), which is invariant under the action of local automorphisms of the $G$-bundles $P$, the so-called  ``gauge group  $\Gau(P)$''.
\end{itemize}
 Given a Lagrangian or   differential equations invariant under a group $G$, we use the expression ``gauging the group $G$'' to indicate the construction of  a gauge theory of $G$, whose Lagrangian or equations reduce to the given ones when the gauge field  is flat.   However, in Physics  the expression ``gauging"  has often a  wider meaning. Basically, it   indicates any process of construction of Lagrangians (or differential equations) that  starts from a given   $G$-invariant Lagrangian and ends up with   a new Lagrangian, which is  invariant under a class of transformations,  locally  identifiable with $G$-valued maps $g: \cU \subset M \to G$ on the space-time $M$. \par
In many classical papers,   it is shown how  General Relativity and other theories of gravity can be obtained by   ``gauging''   the Poincar\`e group $G^\cP$ in  this wider sense or according to a definition of gauge theory different from the above
 (see the classical  \cite{Ki, Ut} or  \cite{Ch, Ne,HMMN, He, SZ} and  the vast bibliography in those papers and books). But not many papers  are concerned with a description of General Relativity  as a gauge theory according to previous  definition (see \cite{Pi, St, IN}). On the other hand, we strongly believe that a sure understanding of this latter description  is a necessary step  if one aims to geometric constructions  of supergravity theories  as gauge theories of  super-extension of the Poincar\`e group  (see \cite{AC} for the classification of these super-extensions  in any dimension and signatures).\par
\smallskip
The purpose of this paper is to give a complete description of how a theory of gravity can be considered as a gauge theory of $G^\cP$. As a by-product we show that, if  a theory of this kind is required to satisfy   the Equivalence Principle, then the  Lagrangian or the differential equations of the theory must be preserved by the gauge transformations {\it plus} an additional  family of transformations, called {\it pseudo-translations}. Finally, we provide a presentation of  such pseudo-translations, which  admits immediate  generalizations to the context of supermanifolds and hence, hopefully, of supergravity theories. \par
\smallskip
Our presentation of gravity as gauge theory is essentially equivalent to the one in \cite{Pi}, but it differs in the following aspect. Any connection form $\o$ on a principal bundle $\pi: P \to M$ induces an horizontal distribution on any associated vector bundle $\pi^E: E \to M$ and hence a differential operator $\n$ of covariant derivation for the sections of $E$. In case $P$ is a $G^\cP$-bundle, there is always an associated bundle $\wt E$, which has the structure of  an affine bundle, endowed with a pseudo-Riemannian metric on the  fibers. Moreover,   there is  a natural one-to-one correspondence between the covariant derivations of sections of $\wt E$  and the connection forms on $P$ (Prop. \ref{2.5}).  Having  applications to supergravity  in mind,  we  represent the gauge fields of a $G^\cP$-theory just as covariant derivations of an affine bundle with metric. By the previous observation, this is  fully equivalent to consider  gauge fields  as in (a) and in \cite{Pi}, but has the advantage  that  deals with objects (covariant derivations) that have  easily defined ``super'' analogues, in contrast with the  ``super'' analogues of principal bundles and connection forms, which require not so straightforward  mathematical  definitions (see e.g. \cite{DM, Al, Sta}). \par
\smallskip
Our first main result can be summarized  as follows. Let $\wt \pi: \wt E \to M$ be an affine bundle over $M$ modeled on a vector bundle $\pi: E \to M$, $(\wt \n, \n)$  a covariant derivation   for the sections of $\wt E$ and $\psi_o: M \to \wt E$ a fixed section of $\wt E$ (see \S 2 for all definitions). If $\psi_o$ satisfies suitable regularity conditions, the  map
 $$L: TM \to E\ ,\qquad L(v) \=  \n_v \psi_o  \in T_x^v \wt E \simeq E_x\qquad\ v \in T_x M\ ,\ x \in M\ , $$ 
 is  a bundle isomorphism that can be used to induce from $\wt E$ a pseudo-Riemannian metric $g$ and a metric covariant derivation  $\n$ on $M$. This construction is used to show that Ê{\it there is a natural correspondence between the pairs $((\wt \n, \n); \psi_o)$, formed by a covariant derivation on $\wt E$  and  regular sections $\psi_o$ of $\wt E$, and the pairs $(g, \n)$, formed by  pseudo-Riemannian metrics and metric covariant derivations on $M$\/}. Since such pairs $(g, \n)$   are the objects usually considered to represent  gravitational fields,  we conclude that any theory of a gravity,  determined  by a Lagrangian on pairs $(g, \n)$, can be interpreted as a gauge theory, provided that:
\begin{itemize}
\item[--] the role of the bundles, whose sections represent the particles coupled with the gauge group, is  played by metric affine bundles  $\wt E$; 
\item[--] the role of the gauge fields is played by   the affine covariant derivations $(\wt \n, \n)$ on $\wt E$; 
\item[--] the theory includes a special field ({\it  Higgs field\/}),  coupled with the gauge field ,  consisting of  a regular section  $\psi_o$ of $\wt E$;
\item[--] the Lagrangian and the constraints, expressed in terms of $(\wt \n, \n)$, are invariant under the family of all gauge transformations of $\wt E$,  explicitly described in \S 2.1.3.
\end{itemize}
We leave to the reader the physical interpretation of the objects involved in this scheme. We need however to stress the following fact: {\it any gauge transformation on a pair  $((\wt \n, \n); \psi_o)$  corresponds just to the identity transformation when expressed in terms of  the pair $(g, \n)$\/}. This is due to the presence of the Higgs field $\psi_o$: It reduces the representation of the  gauge group $\Gau(P)$ to the faithful representation of a smaller class  $\cG_{\psi_o} \subsetneq \Gau(P)$,   which act  on the vielbeins of $g$ by pointwise dependent Lorentz transformations and hence makes  absolutely no change on  $g$ and  $\n$ (see also \cite{St}).
This means that {\it \underline{any} choice for the  Lagrangian and the constraints  on  the pairs $(g, \n)$ corresponds to  a gauge theory on the corresponding pairs $((\wt \n, \n), \psi_o)$\/}. \par
\smallskip
On the other hand, any sensible theory of gravity has to satisfy the so-called {\it Equivalence Principle\/}.  This principle can be briefly stated  saying that the Lagrangian of the theory must be invariant under arbitrary local diffeomorphisms of the space-time (or, equivalent, under arbitrary  local changes of coordinates). A slightly weaker requirement is to ask that the Lagrangian is invariant under 
the diffeomorphisms generated in the  flows of local vector fields. We call it {\it infinitesimal version of Equivalence Principle} and, for simplicity, we consider only this weaker condition. \par
 We show that it corresponds to require that the Lagrangian on the pairs $((\wt \n, \n), \psi_o)$ must be invariant under certain bundle transformations of the metric affine bundle $(\wt E, E)$,  which we call {\it pseudo-translations\/} and of which we determine the explicit action on the gauge field. These pseudo-translations  correspond to the transformations  on $(g, \n)$ determined by  the flows of vector fields on $M$ and  their  expressions in coordinates resembles closely certain gauge transformations. This is probably at the origin of the widespread idea that the transformations,  given by flows of vector fields of the space-time,  are the  outcome  of  a ``gauging process'' of the transformations of the Poincar\`e group.\par
\smallskip
The notion of pseudo-translations bring  to our third result: {\it theories of gravity given by Lagrangians on pairs of the form $(g, \n)$ and satisfying  the infinitesimal version of the Equivalence Principle are in natural correspondence with the gauge theories of the Poincar\`e group $G^\cP$, whose  Lagrangians is not only  invariant  under the gauge group $\Gau(P)$ but also under the family  of all pseudo-translations\/}. In other words, there are {\it two\/} classes of transformations 
to be considered when a gravity theory is presented as a gauge theory \footnote{See also   \cite{IN},  where  the authors adopt  a notion of gauge theory, according to which a covariance  under two (and not one) groups of transformations is required.}. In a future paper, we  will consider the analogous picture in case of  super-extensions of the Poincar\`e group,  providing a  new interpretation of the  ``gauging''   processes of  the super Poincar\`e group   in the construction of theories of supergravity.  \par
\smallskip
The structure of the paper is as follows. In \S 2, we give the first properties of   the covariant derivations on  pseudo-Riemannian affine  bundles  and of  groups of gauge transformations of $G^\cP$-bundles. In \S 3 we  discuss  the theories of gravity as gauge theories of $G^\cP$. In \S 4, we consider the Equivalence Principle,  we introduce the notion  of pseudo-translations 
and prove the stated correspondence  between the infinitesimal Equivalence Principle and a principle of covariance under  pseudo-translations. \par
%
\bigskip
\noindent{\it Notation.} The class of all smooth real functions on a manifold $M$ is denoted by $\gF(M)$, while the class  of vector fields  is denoted by $\gX(M)$. For any bundle $\pi: E \to M$, we denote by $\Sigma(E)$ the family of all global sections  of $E$  and by $\Sigma_{loc}(E)$ the class of local sections defined on open subsets $\cU \subset M$. 
\par
\bigskip

\section{Preliminaries}\label{prel}
\subsection{Gauge transformations of metric affine bundles}

\subsubsection{Affine bundles and their covariant derivations}\label{2.1.1}\hfill\par
\smallskip
\par
\begin{definition} \cite{Go} Let $M$ be a   manifold of dimension $n$. An {\it affine bundle over $M$ modeled on a vector bundle $\pi:E\to M$} is a  fiber bundle $\wt \pi:\wt E\to M$ equipped  with a bundle morphism $+:\wt E\times_{_M}E\to \wt E$
so that, for each $x\in M$, the induced map
\beq\label{morphism} +|_x:\wt E_x\times E_x\to\wt E_x\ ,\qquad \wt e, e \longmapsto \wt e + e\eeq
determines  a structure of affine space on  $\wt E_x $, modeled on the vector space $E_x$. The rank of  $E$ is called {\it rank\/} of the affine bundle.
\end{definition}
In the following, we often denote an affine bundle simply by the pair $(\wt E, E)$. 
Notice that
the bundle map  (\ref{morphism}) induces  the  map $+:  \Sigma(\wt E) \times \Sigma(E) \to \Sigma(\wt E)$,  defined by $(\wt \s +  \s)_x \= \wt \s_x + \s_x$, 
which makes  $\Sigma(\wt E)$  an affine space  modeled on $\Sigma(E)$.
\smallskip
\begin{definition} \label{covder} A {\it covariant derivation on an affine bundle $(\wt E, E)$} is a pair $(\wt\n,\n)$,
where $\n$ is a   linear connection
on $E$ and $\wt\n$ is an operator
$\wt\n:\gX(M)\times \Sigma(\wt E)\longrightarrow\Sigma(E)$
that associates to any vector field $X \in \gX(M)$ and any $\psi \in \Sigma(\wt E)$ a section $\wt \n_X \psi \in \Sigma(E)$ that  satisfies
\begin{itemize}
\item[i)] $\wt\n_{_{fX+gY}}\psi=f\wt\n_{_X}\psi+g\wt\n_{_Y}\psi$;
\item[ii)] $\wt\n_{_X}(\psi+\varphi)=\wt\n_{_X}\psi+\n_{_X}\varphi$
\end{itemize}
for any $\psi\in\Sigma(\wt E)$, $\varphi\in\Sigma(E)$, $X,Y\in\gX(M)$ and
$f,g\in\gF(M)$.\par
\end{definition}
This  definition  is motivated by the following facts.  \par
\medskip
Let us call  {\it affine frame of $\wt E$ at $x$} any  pair $u = (\wt e_0; (e_1, \dots, e_p))$, formed by a point  $\wt e_0$ of  the fiber $\wt E_x$ and a basis  $(e_1, \dots, e_p)$ of the vector space $E_x$.
The family   $A(\wt E)$ of affine frames of $\wt E$ has a natural structure of principal bundle over $M$,  with structure group $\Aff(\bR^p) = \GL_p(\bR) \ltimes \bR^p$, and  $\wt E$ and $E$ are naturally isomorphic to the associated bundles
\beq \wt E \simeq A(\wt E) \times_{\Aff(\bR^p), \wt \rho} \bR^p\ ,\qquad E \simeq A(\wt E) \times_{\Aff(\bR^p),  \rho} \bR^p\eeq
where $\wt \rho$ is the standard representation of $\Aff(\bR^p)$  as group of affine transformations of $\bR^p$  and $\rho$ is the surjective homomorphism  from  $\Aff(\bR^p)$ onto $\GL_p(\bR) = \GL_p(\bR) \ltimes \bR^p/\bR^p$.\par
Any connection form $\omega$ on $A(\wt E)$ determines a parallel transport between  fibers
on $A(\wt E)$, $\wt E$ and $E$.  In particular it defines a pair of  covariant derivations  $\wt \n$ and  $\n$  for the sections of $\wt E$ and $E$, respectively. It can be checked that $\n$ is always  a linear connection, while  $\wt \n$ is an operator that  satisfies (i) and (ii)  of Definition \ref{covder}.  Conversely, for any pair $(\wt \n, \n)$ that satisfies the conditions of Definition \ref{covder}, there exists a  connection form on $A(\wt E)$, whose parallel transport determines $\wt \n$ and $\n$ as associated covariant derivations. In other words,  {\it the covariant derivations $(\wt \n, \n)$ considered  in Definition \ref{covder} coincide with the covariant derivations determined by  the connection forms $\o$ on the principal bundle $A(\wt E)$\/} (see e.g.  \cite{Ta}, \S 2.2).\par
\bigskip
Let $(\wt\n,\n)$ be a covariant derivation on $(\wt E, E)$ and fix a global section  $\psi_o\in\Sigma(\wt E)$ (by well-known properties of fiber bundles,   $\Sigma(\wt E)$ is not empty).  Let
 $L:\gX(M)\to\Sigma(E)$ be  the linear map
 \beq\label{2.3}
L(X)=\wt\n_{_X}\psi_o\ .\eeq   From (i) and (ii) of Definition \ref{covder}, one can check that $L:\gX(M)\to\Sigma(E)$ can be identified with a section $L$ of $T^* M \otimes_M E$.
Given $L$ and $\n$,  the operator $\wt \n$ can be recovered by means of the identity
\beq \label{recovering}
\wt\n_{_X}\psi=L(X)+\n_{_X}(\psi-\psi_o)\ ,\qquad \psi\in\Sigma(\wt E)\ ,\ \ X\in\gX(M)\ .\eeq
Conversely, given  a section $L \in \Sigma(T^* M \otimes_M E)$ and  a linear connection $\n$  on $E$, one can check that the operator $\wt \n: \gX(M) \times \Sigma(\wt E) \to \Sigma(E)$,  defined by  (\ref{recovering}),   satisfies conditions (i) and (ii) of Definition \ref{covder} and  that the correspondence
$$(\wt \n, \n) \longrightarrow (L \= \wt \n \psi_o, \n)$$
is one to one. For this reason, once  $\psi_o$ is given, we will  often indicate a  covariant derivation $(\wt\n,\n)$ by the corresponding pair  $(L,\n)$.\par
\medskip
\subsubsection{Pseudo-Riemannian  affine bundles and metric connections}\hfill\par
\smallskip
We now  introduce the concept of pseudo-Riemannian metrics on affine bundles. Let $\pi: Q \to M$ be a fiber bundle
and denote by $T^v Q \subset TQ$  the vertical distribution.
  We call {\it vertical tensor field of type $(r,s)$\/} on $Q$ any section of the bundle 
$ \pi: \bigotimes^r T^v Q \otimes_Q \bigotimes^s T^{v*} Q \longrightarrow Q$.
A vertical tensor field $\wt \a$ of type $(r,s)$  of an  affine  bundle $\pi: \wt E \to M$ modeled on $E$ will be called {\it affine} if, for any $x \in M$,  there exists a tensor $ \a_x \in \bigotimes^r E_x \otimes \bigotimes^s E_x$ so that for any $ u \in E_x$  the tensor $\wt \a_u$ is identifiable with  $\a_x $   under the natural isomorphism $T^v_u \wt E \simeq E_x$.  Therefore any affine tensor field $\wt \a$ on $\wt E$ is uniquely determined by   a corresponding section $\a$ of the tensor bundle
$ \pi: \bigotimes^r E \otimes_M \bigotimes^s E^* \longrightarrow M$.
\begin{definition} A {\it pseudo-Riemannian metric\/} on an affine bundle $(\wt E, E)$  is an affine vertical tensor field $\wt g$ of type $(0,2)$,  determined by   a section $g$ of $\bigotimes^2 E^*$, so that   $g_x$ is an inner product on $E_x$ of constant signature for any $x \in M$. The triple $(\wt E, E, g)$ is called {\it pseudo-Riemannian affine bundle modeled on $E$\/}.\par
We also call {\it pseudo-Riemannian metric on a vector bundle $E$\/}ÿ any  vertical tensor field of type $(0,2)$ on $E$ that satisfies the above conditions, taking  $E$ as  affine bundle modeled on itself.
\end{definition}
\par
Given a pseudo-Riemannian affine bundle $(\wt E, E, g)$ of signature $(p,q)$,  any affine frames $u = (\wt e_0, (e_1, \dots, e_p))$ with  $(e_1, \dots, e_p)$ orthonormal w.r.t. the inner product $g_x$, $x = \pi(u)$, is called {\it orthonormal frame bundle\/} and the class $O_g(\wt E)$ of all orthonormal affine frames is a $O_{p,q}\ltimes \bR^{p+q}$-reduction of $\Aff(\wt E)$. Conversely, for any  $O_{p,q}\ltimes \bR^{p+q}$-reduction  $P \subset \Aff(\wt E)$, there exists a unique pseudo-Riemannian metric $\wt g$ on the associated bundle $\wt E = O_g(\wt E) \times_{O_{p,q} \ltimes \bR^{p+q}, \wt \r} \bR^{p+q}$ for which $P$ is the corresponding orthonormal affine bundle $P = O_g(M)$.\par
\medskip
\begin{definition} Let $(\wt E, E, g)$ be a pseudo-Riemannian affine bundle and $(\wt \n, \n)$ an affine covariant derivation of $(\wt E, E)$. We say that $(\wt \n, \n)$  is {\it metric\/} if the linear connection $\n$ is so that
$\n_X g = 0$ for any $X \in \gX(M)$.
\end{definition}
The following proposition can be obtained with standard arguments on reductions. For a detailed proof, see e.g.  \cite{Ta}, \S 3.2.2.\par
\medskip
\begin{prop}\label{2.5} Any connection form $\o$ on the orthonormal affine frame bundle $O_g(\wt E)$ of a
 pseudo-Riemannian affine bundle $(\wt E, E, g)$  induces a pair  covariant derivation $(\wt \n, \n)$ on the associated bundles $\wt E = O_g(\wt E) \times_{O_{p,q}\ltimes \bR^{p+q}, \wt \r} \bR^{p+q}$  and
 $E = O_g(\wt E) \times_{O_{p,q}\ltimes \bR^{p+q}, \r} \bR^{p+q}$, which is a metric affine covariant derivation. Conversely, any  metric affine covariant derivation $(\wt \n, \n)$
is uniquely determined by a connection form $\o$ on $O_g(\wt E)$.
\end{prop}
\bigskip

\subsubsection{Affine gauge transformations}\label{2.1.3}\hfill\par
\bigskip
We recall that a gauge transformation of a principal bundle $\pi:P\to M$ is a bundle
automorphism $f:P\to P$ which induces the identity on $M$. The gauge transformations of
the associated bundles are the bundle automorphisms induced by the gauge transformations
of $P$. When $P$ is the orthonormal frame bundle $P=O_g(\wt E)$ of a pseudo-Riemannian affine bundle $(\wt E,E,g)$, one can check that
the gauge transformations of $\wt E$ and $E$ coincide with the maps given in the following definition.\par
\begin{definition}
Let $(\wt E,E,g)$ be a pseudo-Riemannian affine bundle. A {\it gauge transformation of $E$} (resp. $\wt E$) is any bundle automorphism $f$ of $E$ (resp. $\wt E$) which preserve the fibers and,  for any $x\in M$,  the map
$\left.f\right|_{E_x}$ (resp. $\left.f\right|_{\wt E_x}$) is an isometry of $(E_x,g_x)$ (resp. $(\wt E_x,g_x)$).\par
We call {\it infinitesimal gauge transformation of $E$} (resp. {\it of $\wt E$}) any vector field $V\in\gX(E)$ (resp. $\gX(\wt E)$) whose
flow $\Phi^V_t$ is a one-parameter group of gauge transformations.
\end{definition}

We denote by $\Gau(\wt E)$ and $\Gau(E)$ the groups of gauge transformations of $\wt E$ and $E$, respectively. We remark that there exists a natural surjective homomorphism
$\alpha:\Gau(\wt E)\to\Gau(E)$ defined as follows. If $f\in\Gau(\wt E)$ we have that $f|_{\wt E_x}$ is an isometry of $(\wt E_x,g_x)$ and hence it is an affine transformation for any $x\in M$. Let $h_x:E_x\to E_x$ be the linear map associated to $f|_{\wt E_x}$ so that
\beq\label{alpha}
f(\wt e')-f(\wt e)=h_x(\wt e'-\wt e)
\eeq
for any $\wt e,\wt e'\in\wt E_x$. The map $h=\alpha(f)\in\Gau(E)$ is the gauge transformation,  which induces on each fiber $E_x$ the linear map $h_x$ just described.
Clearly, $f\in\ker\alpha$ if and only if $f|_{\wt E_x}$ is a translation of the affine space $\wt E_x$ for any $x\in M$, i.e. if and only if it is a map of the form
\beq
\wt\xi:\wt E\to\wt E, \qquad \wt\xi(\wt e)\=\wt e+\xi|_x, \quad x=\wt\pi(\wt e).
\eeq
for some section $\xi\in\Sigma(E)$. We call any map of this kind {\it point-depending translation ({\rm or simply} translation) determined by $\xi$} and we denote by $\Gau^T(\wt E)$
the normal subgroup of $\Gau(\wt E)$ of all translations.
Summing up, we have the following exact sequence:
\beq\label{sequence}
0 \longrightarrow \Gau^T(\wt E) \overset{\imath} \longrightarrow \Gau(\wt E) \overset{\alpha} \longrightarrow \Gau(E) \longrightarrow 1,
\eeq
which is the analogue of the corresponding sequence for $\Aff(\bR^p)$.\par

Now, fix a section $\psi_o\in\Sigma(\wt E)$. A gauge transformation $h:E\to E$ induces the following gauge transformation on $\wt E$
\beq
\wt h^{\psi_o}:\wt E\to\wt E, \qquad \wt h^{\psi_o}(\wt e)\=\psi_o(x)+h(\wt e-\psi_o(x)), \quad x=\wt\pi(\wt e).
\eeq
This map will be called {\it point-depending rotation ({\rm or simply} rotation) around $\psi_o$ determined by $h$}.

\begin{lem}\label{lemmino}\hfill\par
\begin{itemize}
\item[1)] A gauge transformation $f:\wt E\to\wt E$ is a rotation around $\psi_o$ if and only if $f(\psi_o(x))=\psi_o(x)$ for any $x\in M$.
\item[2)] Any gauge transformation $f$ of $\wt E$ can be uniquely written as
\beq
f=\wt\xi\circ\wt h,
\eeq
where $\wt\xi$ is the translation by the section $\xi(x)\=f(\psi_o(x))-\psi_o(x)$, while $\wt h$ is the rotation around $\psi_o$ defined by $\wt h\=\wt\xi^{-1}\circ f$.
\item[3)]
The homomorphism $\beta:\Gau(E)\to\Gau(\wt E)$ defined by $\beta(h)\=\wt h^{\psi_o}$ makes (\ref{sequence}) a splitting exact sequence.
\end{itemize}
\end{lem}
\begin{pf}
(1) The necessity is immediate. Conversely, assume that $f$ fixes $\psi_o$ and let $h=\alpha(f):E\to E$ be the gauge transformation defined by (\ref{alpha}).
From construction it follows that
\[
f(\wt e)=f(\psi_o(x))+h(\wt e-\psi_o(x))=\psi_o(x)+h(\wt e-\psi_o(x)), \quad x=\wt\pi(\wt e),
\]
and hence that $f$ coincides with the rotation around $\psi_o$ by $h$.\par
\noindent (2) We only need to check that $\wt h$ fixes $\psi_o$. Since $\wt\xi^{-1}$ is the translation by the section $-\xi=\psi_o-f\circ\psi_o$, the conclusion follows immediately.\par
\noindent (3) It is immediate to verify that $\alpha\circ\beta=\mathrm{Id}$.
\end{pf}

Recall that any vertical tangent subspace $T^v_{\wt e}\wt E$ is naturally identified with the fiber $E_x$, $x=\wt\pi(\wt e)$, of $E$. By this identification, any section $\xi\in\Sigma(E)$ corresponds to a vertical vector field $V^\xi$ on $\wt E$ whose flow $\Phi^{V^\xi}_t$ consists of the translations by the sections $t\cdot\xi\in\Sigma(E)$, $t\in\bR$. In other words,  $V^\xi$ is an infinitesimal gauge transformation and  we call it {\it infinitesimal translation}.\par

\subsubsection{Actions of gauge transformations on affine covariant derivations}
\label{gauandcov}
\begin{definition}\label{infdef}
Let $(\wt\n,\n)$ be an a metric affine covariant derivation on $(\wt E,E,g)$ and $f\in\Gau(\wt E)$. Let also $h=\alpha(f)$ the corresponding gauge transformation on $E$.
We call {\it deformation of $(\wt\n,\n)$ by $f$} the covariant derivation $(f_*\wt\n,f_*\n)$ defined by
\beq\label{def_nabla}
f_*\wt\n\=h\circ\wt\n\circ f^{-1}, \qquad f_*\n\=h\circ\n\circ h^{-1}.
\eeq
If $V\in\gX(\wt E)$ is an infinitesimal gauge transformation, we call {\it infinitesimal deformation of $(\wt\n,\n)$ by $V$} the pair of operators $(V(\wt\n),V(\n))$
\[
V(\wt\n):\gX(M)\times\Sigma(\wt E)\to \Sigma(E), \qquad V(\n):\gX(M)\times\Sigma(E)\to \Sigma(E)
\]
respectively defined by
\[
V(\wt\n)_{_X}\wt e\=\left.\frac{d}{dt} \left(\Phi^V_t{}_*\wt \n\right)_{_X} \wt e\right|_{t=0}, \qquad V(\n)_{_X} e\=\left.\frac{d}{dt} \left(\Phi^V_t{}_*\n\right)_{_X} e\right|_{t=0}
\]
for any $X\in\gX(M)$, $\wt e\in\Sigma(\wt E)$ and $e\in\Sigma(E)$.
\end{definition}

\begin{rem}
By Proposition \ref{2.5}, any metric affine covariant derivation on $(\wt E,E)$ corresponds to a unique connection form $\omega$ and associated horizontal distribution $\cH=\ker\omega$ on $O_g(\wt E)$. On the other hand, by the remarks in \S \ref{2.1.3}, $f\in\Gau(\wt E)$ corresponds to a unique gauge transformation $\hat f:O_g(\wt E)\to O_g(\wt E)$.  One can check that $(f_*\wt\n,f_*\n)$ coincides with the affine covariant derivation associated with the push-forward of $\cH$ by $\hat f$, which is a horizontal distribution associated with the connection form $\hat f^*\omega$ (see e.g. \cite{Ta}). This fact motivated the  previous Definition \ref{infdef}.
\end{rem}

We recall that, by the remarks in \S \ref{2.1.1},
given a fixed section $\psi_o\in\Sigma(\wt E)$, the covariant derivation $(\wt\n,\n)$ can be identified with the pair $(L,\n)$, with $L$
defined in (\ref{2.3}). A straightforward computation shows that the pair $(f_*L,f_*\n)$, corresponding to the deformation $(f_*\wt\n,f_*\n)$ by
$f\in\Gau(\wt E)$, is given by
\beq\label{def_L}
(f_*L)(X)=h\left(L(X)+\n_{_X}\left(f^{-1}(\psi_o)-\psi_o\right)\right)
\eeq
and of course by the second formula in (\ref{def_nabla}). It follow immediately that,
if $f$ is a translation by $\xi\in\Sigma(E)$, the deformation $(f_*L,f_*\n)$ is
\beq\label{def_tran}
f_*L=L-\n\xi, \qquad f_*\n=\n,
\eeq
while, if $f$ is a rotation around $\psi_o$ by $h\in\Gau(E)$, $(f_*L,f_*\n)$ is
\beq\label{def_rot}
f_*L=h\circ L, \qquad f_*\n=h\circ\n\circ h^{-1}.
\eeq

\bigskip
\section{Gravity theories as theories of metric affine connections}
\subsection{Gravitational fields and metric affine connections}\hfil\par
\setcounter{equation}{0}
\medskip
\label{classicalgravityfields}
Let $M$ be a manifold of dimension $n$. We call {\it gravitational field on $M$\/} any pair $(g, \n)$ formed by a pseudo-Riemannian metric $g$ of signature $(p,q)$ and a metric covariant derivation $\n$, i.e. so that $\n g = 0$.  This terminology is motivated by the fact that, in General Relativity, gravity is represented by a 4-dimensional space-time $M$ and a pair $(g, \n)$, where $g$ is a  pseudo-Riemannian metric  of signature $(1,3)$ and $\n$ is the Levi-Civita connection of $g$ (i.e. metric and torsion free)  so that  the well-known Einstein equations are satisfied.  In other words, we may say   that {\it in General Relativity the gravity is represented by a gravitational field $(g, \n)$  of signature $(1,3)$ satisfying the conditions 
\beq \label{GR} T = 0, \qquad Ric - \frac{s}{2} g =  \mathcal T^g\eeq
where $T$, $Ric$ and $s$ are  the torsion, the Ricci tensors and the scalar curvature of $\n$, respectively,  and  $\mathcal T^g$
the stress-energy tensor determined by other physical fields\/}. 
It is therefore natural to consider the generalizations of General Relativity   as theories of gravitational fields $(g, \n)$, subjected to systems of equations that are extensions or modifications of  (\ref{GR}). \par 
\medskip
In this section, we want to show how any gravitational field $(g, \n)$ can be naturally associated with a covariant derivation on a suitable metric affine bundle, i.e. to the associated bundle of a principal bundle with structure group given by the Poincar\`e group $G = O_{p,q} \ltimes \bR^n$ (see also \cite{Mc, PD}).\par
\bigskip
Let  $(\wt E,E,g_o)$ be  a metric affine bundle over $M$ of rank $n$ and $\psi_o$ a section of $\wt E$.
 A metric affine covariant derivation $(L,\n)$
will be  called {\it regular w.r.t. $\psi_o$\/} if the map $L_x:T_xM\to E_x$ is a linear isomorphism for any $x\in M$.  Clearly, if a quadruple $(\wt E,E,g_o, \psi_o)$ admits a regular covariant derivation, then there exists an affine bundle isomorphism between $(\wt E, E)$  and  $(TM, TM)$ which maps  $\psi_o$ into the zero section $\wt 0$ of $TM$ (here we consider $TM$ as an affine bundle modeled on itself). \par
\smallskip
From now on we will always assume that $(\wt E,E,g_o)$ admits some regular covariant derivation and hence that $(\wt E,E, g_o, \psi_o)\simeq(TM,TM, g_o, \wt 0)$ for some pseudo-Riemannian metric $g_o$ on $M$. However,  since  the (non canonical) identification $\wt E \simeq TM$, $E \simeq TM$, etc.  is not relevant and it might even  cause confusion in certain arguments,   we will avoid it in all what follows.\par
\bigskip
As observed in \S \ref{2.1.1},  any given  $\psi_o \in \Sigma(\wt E)$ brings  to the identification of covariant derivations $(\wt \n, \n)$  with corresponding pairs $(L,\n)$. On the other hand, any such pair  
associated with a regular covariant derivation determines a gravitational field $(g, \n^L)$ on $M$ as follows:
\beq\label{corr}
g(X,Y)\=g_o(L(X),L(Y)), \qquad \n^L_{_X}Y\=L^{-1}(\n_{_X}(L(Y))).
\eeq
Notice that $\n^L$ is metric w.r.t.  $g$, because   $\n$ metric w.r.t.  $g_o$ and hence $\n^L_{_X}g=0$ for any  $X \in \gX(M)$. \par
Using (\ref{def_rot}), one can check that two regular covariant derivations $(L,\n)$, $(L',\n')$
determine the same gravitational field $(g,\n^L)$ if and only if they differ by a rotation around $\psi_o$.  \par
Now, let us  consider the following notation: 
\begin{itemize}
\item[-]  $\Conn(\wt E)$ denotes the class of all metric affine covariant derivations of $(\wt E, E, g_o)$; 
\item[-]  $\Conn(\wt E)^{\psi_o\reg} \subset \Conn(\wt E)$ is the subclass of regular ones w.r.t. $\psi_o$
and $\Conn(\wt E)^{\reg} = \bigcup_{\psi_o\in \Sigma(\wt E} \Conn(\wt E)^{\psi_o\reg}$; 
\item[-] $\cG_{\psi_o} \= \Gau_{\psi_o}(\wt E) \subset \Gau(\wt E)$ is the isotropy subgroup of $\Gau(\wt E)$  at $\psi_o$ (i.e.  the group of rotations around $\psi_o$); 
\item[-] $\Grav_{p,q}(M)$ is the class of  all gravitational fields  $(g, \n)$ on $M$  with $g$ of  signature $(p,q)$ (the same of $g_o$) . 
\end{itemize}
Then, the  correspondence  described in (\ref{corr}) determines a map 
\beq \label{projection} \imath_{\psi_o}: \Conn(\wt E)^{\psi_o\reg} \longrightarrow \Grav_{p,q}(M)\eeq
which induces an injection from  the space of  $\cG_{\psi_o}$-orbits  $\Conn(\wt E)^{\psi_o\reg} /\cG_{\psi_o}$  into the family of gravitational field.\par
\smallskip
The  map $\imath_{\psi_o}$ is indeed a projection. In fact, by standard facts on inner products, any pseudo-Riemannian metric $g$ of signature $(p,q)$  is of the form (\ref{corr}) for some suitable   tensor field $L$ of type $(1,1)$. Moreover, if $\n'$ is a covariant derivation on $M$ which is metric for $g$, then the covariant derivation on $E$ defined by $\n = L \circ \n' \circ L^{-1}$ is metric for $g_o$ and the pair $(L, \n)$ is mapped onto $(g, \n')$ via  (\ref{corr}). This proves the surjectivity.
Summing up,  we proved the following. 
\begin{theo}\label{firsttheorem} For any given $\psi_o \in \Sigma(\wt E)$, 
the map $\imath_{\psi_o}$  defined  in (\ref{corr}) induces a  one to one correspondence between the orbit space $\Conn(\wt E)^{\psi_o\reg}/\cG_{\psi_o} $ and  gravitational fields  $(g, \n)$ in  $\Grav_{p,q}(M)$.
\end{theo}
\begin{coord}
Let $(x^1,\ldots,x^n):\cU\subset M\to \bR^n$ be a system of coordinates on $M$ and fix a collection $(e_1^o,\ldots,e_n^o)$ of sections $e_i^o\in\Sigma(E)$ so that
$(e^o_i|_x)$ is an orthonormal basis of $(E_x,g_o)$. Any metric affine covariant derivation $(L,\n)$ is of  the form
\[
L=\theta^i_\mu e^o_i \otimes dx^\mu, \qquad \n_{\frac{\partial}{\partial x^\mu}}(\varphi^i e^o_i)=\left(\frac{\partial\varphi^i}{\partial x^\mu}+\Gam i \mu j \varphi^j\right) e^o_i,
\]
where $\Gam i \mu j$ are the components of the derivatives $\n_{\frac{\partial}{\partial x^\mu}}e^o_j=\Gam i \mu j e^o_i$.
The vector fields 
\[
e_i=e_i^\mu\frac{\partial}{\partial x^\mu}\=L^{-1}(e^o_i)\in\gX(M)
\]
constitute an orthonormal frame field (or {\it vielbein}) for the metric $g=g_o(L(\cdot),L(\cdot))$, while $\n^L$ is of the form
\[
\n^L_{\frac{\partial}{\partial x^\mu}}(X^i e_i)=\left(\frac{\partial X^i}{\partial x^\mu}+\Gam i \mu j X^j\right) e_i
\]
Since $e_i^\mu \theta_\mu^j=\delta^j_i$, it is clear that the application
$L$ can be completely recovered from the vielbein $(e_i)$. Moreover,
\beq\label{gamma}
0=g(\n^L_{\frac{\partial}{\partial x^\mu}} e_i,e_j)+g(e_i,\n^L_{\frac{\partial}{\partial x^\mu}} e_j)=\Gam j \mu i+\Gam i \mu j
\eeq
and one can easily check that any set of functions $\Gam j \mu i$ that satisfies (\ref{gamma}) determines uniquely a metric covariant derivation $\n'$ on $M$ and hence a
metric covariant derivation $\n=L\circ\n'\circ L^{-1}$ on $E$. Therefore, the class of covariant derivations $(L,\n)$ can be locally identified with the pairs $((e_i),(\Gam j \mu i))$ formed by a vielbein $(e_i)$ and functions $\Gam j \mu i$ with
$\Gam j \mu i=-\Gam i \mu j$.\par
\medskip
Let us now write the formulae that express the action of the gauge transformations in terms of the pairs $((e_i),(\Gam j \mu i))$. Assume that $f=\wt\xi\in\Gau^T(\wt E)$
is a translation by $\xi=\xi^i e^o_i\in\Sigma(E)$. Denoting as before by $\theta^i_\mu$
the functions which give the components of $L$ and are hence defined in terms of the vielbein $(e_i)$ by the relations $e^\mu_i\theta^j_\mu=\delta^j_i$, one can
immediately obtain that
\beq
((e_i),(\Gam j \mu i))\overset{\wt\xi}\longmapsto((e'_i),({\Gamma'}^{\ j}_{\!\!\mu i}))
\eeq
where ${\Gamma'}^{\ j}_{\!\!\mu i}=\Gam j \mu i$ and $e'_i={e'}^\mu_i\frac{\partial}{\partial x^\mu}$ is defined by the following equations
\beq
{e'}^\mu_i\left(\theta^j_\mu-\frac{\partial\xi^j}{\partial x^\mu}-\Gam j \mu \ell \xi^\ell\right)=\delta^j_i.
\eeq
In case $f=\wt h\in\Gau(\wt E)$ is a rotation around $\psi_o$ by $h\in\Gau(E)$, $h(e^o_i)=h^j_i e^o_j$, we have that
\beq\label{coord_rot}
((e_i),(\Gam j \mu i))\overset{\wt h}\longmapsto\left(\left(\left(h^{-1}\right)^j_i e_j\right),\left(
h^j_\ell \Gam \ell \mu k \left(h^{-1}\right)^k_i+h^j_\ell\frac{\partial\left(h^{-1}\right)^\ell_i}{\partial x^\mu}
\right)\right).
\eeq
Notice that Theorem \ref{firsttheorem} implies that locally any gravitational fields $(g,\n)$ can be identified with a pair of the form $((e_i),(\Gam j \mu i))$
uniquely determined up to a transformation (\ref{coord_rot}).
\end{coord}
\bigskip
\subsection{Theories of gravity as  gauge theories}\label{tggt}\hfill\par
\medskip
Consider the actions of  the elements $f \in \Gau(\wt E)$ on the section $\psi_o$ and on the corresponding map $\imath_{\psi_o}$.   We claim  that
\beq \label{invariance}
 \imath_{f (\psi_o)} \circ f_* = \imath_{\psi_o}\eeq
 for any $f \in \Gau(\wt E)$.
To prove this, by   Lemma \ref{lemmino} and previous remarks, we may assume with  no loss of generality that   $f = \wt \xi$ is  a translation determined by a section $\xi$ of $E$. 
Then
$$f_*(L, \n) = (L - \n \xi, \n)$$ 
$$\imath_{f (\psi_o)} (f_*(L, \n)) = (g_o(\wh L(\cdot), \wh L(\cdot)), \n)\ ,\ \text{where}\ \   \wh L \= (f_*\wt \n) f(\psi_o)\ .$$
Since $ (f_*\wt \n) f(\psi_o) = \wt \n (\wt \xi^{-1}( \wt \xi(\psi_o))) = \wt \n \psi_o$, we have  that $\wh L = L$ and (\ref{invariance}) follows. Due to  (\ref{invariance}) and Theorem \ref{firsttheorem}, if we set 
$\O \= \bigcup_{\psi_o \in \Sigma(\wt E)}\Conn(\wt E)^{\psi_oreg} \times \{\psi_o\}$, 
the map 
\beq \imath:  \Omega\to \Grav_{p,q}(M)\ ,\qquad \quad \label{projection1}  \imath((\wt \n, \n); \psi_o) \= \imath_{\psi_o}(L, \n) = (g, \n^L)\eeq
 induces a one-to-one correspondence between $\O/\Gau(\wt E)$ and $\Grav_{p,q}(M)$. \par
This identification $\Grav_{p,q}(M)\simeq  \O/\Gau(\wt E)$ is not in contrast with the identification $\Grav_{p,q}(M)\simeq \Conn(\wt E)^{\psi_o\reg}/\cG_{\psi_o}$ given in  Theorem \ref{firsttheorem}. In fact, the space $\O$ is union of the $\Gau^T(\wt E)$-orbits of  $ \Conn(\wt E)^{\psi_o\reg} \simeq \Conn(\wt E)^{\psi_o\reg} \times\{\psi_o\}$  and hence 
$\Grav_{p,q}(M)\simeq \Conn(\wt E)^{\psi_o\reg}/\cG_{\psi_o} \simeq \left(\frac{\O}{\Gau^T(\wt E)}\right)/\cG_{\psi_o} =\frac{ \O}{\Gau(\wt E)}$.\par
For practical purposes, the identification $\Grav_{p,q}(M)\simeq  \Conn(\wt E)^{\psi_o\reg}/\cG_{\psi_o}$ is more efficient and it is the only one we use in the following. On the other hand, this last identification $\Grav_{p,q}(M)\simeq \O/\Gau(\wt E)$ allows to state   that  
 {\it the theories on gravity fields $(g, \n)$ are in natural correspondence with  theories on  the triples $(\wt \n, \n,  \psi_o) \in \O$ that are  invariant  under  the full gauge group $\Gau(\wt E)$, i.e.  of the gauge group of the  $O_{p.q} \ltimes \bR^n$-bundle $P = O_g(\wt E)$\/} \footnote{On invariance under the transformations in $\Gau(\wt E)$,  see also \cite{Ta} and \cite{St}, \S 1.}.  Moreover, it must be stressed  that, via the map $\imath$, the action of $\Gau(\wt E)$ on $\O$ corresponds to the trivial action on   $\Grav_{p,q}(M)$.  Hence,  the presentation of the theories on gravity fields as gauge theories of $\Gau(\wt E)$ does not carry  any practical advantage for analyzing  the dynamics of  gravity fields $(g, \n)$. The main application we have in mind is  to  provide a solid scheme of geometric construction for gauge theories with gauge group  of super-extensions of the Poincar\`e group $O_{p.q} \ltimes \bR^n$, i.e. of theories of supergravity. \par
\begin{rem} \label{remark1} Since $\Gau(\wt E)$ acts transitively on the  sections of $\wt E$,  no constraint on $\psi_o$ might occur if one look for gauge invariant   equations on  $\O$. 
On the other hand, since $\Gau(\wt E)$ acts trivially on   $\Grav_{p,q}(M)$,  there is  no effect  if we break the gauge  invariance  by  considering $\psi_o$ as fixed (see also \cite{St}, end of \S 1). So,   with no loss of generality,  we may  state that  {\it the theories on gravity fields $(g, \n)$ are in natural correspondence with the  theories on regular metric affine connections $(\wt \n, \n) = (L, \n)$, which are  invariant under the reduced gauge group $\cG_{\psi_o} = \Gau_{\psi_o}(\wt E)$\/}. 
\end{rem}
\begin{rem} \label{remark2}  Even if the subgroup of translations $\Gau^T(\wt E) \subset \Gau(\wt E)$  does act (locally) on the 
class of regular metric affine covariant derivations $\Conn(\wt E)^{\psi_o reg}$, the map $\imath$ cannot be used to induce  any corresponding action (not even the trivial one) on the space of gravity fields $\Grav_{p,q}(M) \simeq \Conn(\wt E)^{\psi_o r}/\cG_{\psi_o}$.  In fact, {\it  the isotropy gauge group
 $\cG_{\psi_o} = \Gau_{\psi_o}(\wt E)$ is not normalized by the action of  $\Gau^T(\wt E)$\/} and hence there is no induced action of $\Gau^T(\wt E) $ on the quotient $\Conn(\wt E)^{\psi_o r}/\cG_{\psi_o}$. To check this directly, consider two covariant derivations  $(L,\n)$ and $(L',\n')$ in the same $\cG$-orbit (i.e. $L'=h\circ L$ and $\n'=h\circ\n\circ h^{-1}$ for some rotation $h \in \cG$)  and let $\wt\xi \in \Gau^T(\wt E)$   determined  by $\xi \in \Sigma(E)$. Then
$
(\wt\xi_*L,\wt\xi_*\n)=(L-\n\xi,\n)$ and  $(\wt\xi_*L',\wt\xi_*\n')=(h\circ (L-\n\eta),h\circ \n\circ h^{-1})$,
with $ \eta= L^{-1}(\xi)$. 
It follows  that,  in general,  $(\wt\xi_*L,\wt\xi_*\n)$ and $(\wt\xi_*L',\wt\xi_*\n')$ are not in the same $\cG$-orbit. 
\end{rem}
\bigskip
\section{Theories of gravity as\\ gauge theories satisfying the Equivalence Principle}\hfill\par
\setcounter{equation}{0}
\label{EquivalenceP}
The Equivalence Principle of General Relativity (i.e.  covariance under changes of coordinates)  requires that {\it
all constraints and equations  must be covariant under any local diffeomorphism\/}, that is  the class of their solutions  is invariant under the action of local diffeomorphisms. This  corresponds to an invariance  property on the corresponding  gauge theory on  $\Conn(\wt E)^{\psi_o reg}$  {\it distinct\/} from the invariance w.r.t.  to  $\cG_{\psi_o}$,  described in the previous section. 
Therefore, possible  generalizations of General Relativity that satisfy  the Equivalence Principle  must be searched {\it amongst $\cG_{\psi_o}$-invariant theories in $\Conn(\wt E)^{\psi_o reg}$  that are invariant under an additional pseudogroup of local transformations\/}, namely under a  pseudogroup acting on $\Conn(\wt E)^{\psi_o reg}$ in a way that corresponds  to the action of the local diffeomorphisms on  $\Grav_{p,q}(M)$.  In the next two sections  we determine the infinitesimal transformations of such pseudogroup.\par
\bigskip
\subsection{Pseudo-translations and the Equivalence Principle for a theory in $\Conn(\wt E)^{\psi_o reg}$ } \hfill\par
\subsubsection{Torsion and curvature of  metric affine covariant derivation. Parameterizations  by the torsion}\hfill\par
\medskip
\begin{definition} \label{deftor} We call {\it torsion\/} of a regular metric affine covariant derivation $(L, \n)$ the 
section of $\Sigma(\Lambda^2 E^* \otimes E)$ defined by 
$$T(s,s') \= \n_{L^{-1}(s)} s' -  \n_{L^{-1}(s')} s - [L^{-1}(s), L^{-1}(s')]$$
for any $s, s' \in \Sigma(E)$.  We call {\it curvature\/} of  $(L, \n)$ the section of $\Sigma(\Lambda^2 E^* \otimes \so(E, \g_o))$ defined by
$$R(s,s')(s'')  \= \n_{L^{-1}(s)} \n_{L^{-1}(s')} s''  -   \n_{L^{-1}(s')} \n_{L^{-1}(s)} s'' - 
 \n_{[L^{-1}(s), L^{-1}(s')]} s''$$
 for any $s, s',s'' \in \Sigma(E)$.
\end{definition}
Notice  that the tensor fields $T^L = L^* T$ and $R^L = L^* R$ on $M$ coincide with the torsion and curvature of the connection $\n^L = L^{-1} \n \circ L$ associated with $(L, \n)$. This motivates our terminology.\par
\medskip
We want now to show that the torsions can be used to completely parameterize the space of metric affine covariant derivations, in full analogy with the parameterization by torsions of the metric connections on  pseudo-Riemannian manifolds. First of all, 
fix a  metric covariant derivation  $\n^o$  on the vector bundle $(E, g_o)$ and for any  $(L, \n) \in \Conn^{\psi_o reg}(\wt E)$ let 
$$\d \n: \Sigma(E) \times \Sigma(E) \longrightarrow \Sigma(E)$$
\beq \d \n (s,s') \= \n_{L^{-1}(s)} s' - \n^o_{L^{-1}(s)}(s')\ , \qquad s, s' \in \Sigma(\wt E)\ .\eeq
By construction and definitions,   for any $\l,\mu \in \gF(M)$ and $s, s', s''\in \Sigma(E)$ we have that
$\d \n (\l s + \mu s', s'') = \l \d \n (s, s'') + \mu \d \n(s', s'')$ and 
$\d \n (s, \l  s' + \mu s'') = \l \d \n (s, s') + \mu \d \n (s, s')$. 
This means that  $\d \n $ can be  uniquely represented as a vertical tensor field with values in $E^* \otimes_M E^* \otimes_M E$. Being $\n^o$ and $\n$ both metric w.r.t. $g_o$, from equalities 
$$  g_o(\n^o_{L^{-1}(s)}s' - \n_{L^{-1}(s)}s', s'') =$$
$$ =  L^{-1}(s)g_o(s', s'') - 
g_o(s', \n^o_{L^{-1}(s)} s'') - L^{-1}(s)g_o(s', s'') + g_o(s', \n_{L^{-1}(s)} s'') =$$
$$ = g_o(s', \n^o_{L^{-1}(s)}s'' - \n_{L^{-1}(s)}s'')\ ,$$
we conclude that
\beq \label{skew} g_o(\d \n (s, s'), s'') +  g_o(s', \d \n (s,  s'')) = 0\ , \eeq
meaning  that $\d \n$ is indeed a section of $E^* \otimes \so_{g_o}(E, g_o)$ (here 
 $\pi:  \so(E, g_o) \to M$ is the vector bundle of vertical tensor fields in $E^* \otimes E$, skew-symmetric w.r.t. $g_o$, i.e.   with fibers equal to $\so(E_x, g_o)$).
\par
\medskip
Conversely, given a pair $(L, \d \n)$, consisting of vertical  tensor fields  in  $ (T^*M \otimes_M E) + (E^* \otimes\so(E, g_o))$  with  $L_x: T^*_x M \to E_x$ invertible for any $x \in M$,  we may consider the 
pair $(L, \n)$ with $\n$   defined by 
\beq\label{affinestructure}\n_X s \= \n^o_X s + \d \n (L(X),  s). \eeq
From previous remarks, $(L, \n)$ is a metric affine covariant derivation in $\Conn^{\psi_o reg}(\wt E)$ and 
\eqref{affinestructure} allows to $\Conn^{\psi_o reg}(\wt E) \simeq \Sigma\left(T^*M \otimes_M E\right)\times\Sigma\left( E^* \otimes\so(E, g_o)\right)$. 
In particular,  we may conclude that $\Conn^{\psi_o}(\wt E)$ has a structure of Frechet space with tangent spaces isomorphic to the space of 
sections $\Sigma(T^*M \otimes_M E \oplus_M E^* \otimes_M\so(E, g_o))$. \par
Let us now  consider  the skew-symmetrizing map $\partial: \Sigma(E^* \otimes\so(E, g_o)) \longrightarrow \Sigma(\Lambda^2 E^* \otimes E)$, called {\it Spencer operator\/}, defined by
\beq \partial H(s,s') =  H(s, s') -  H(s', s) \ .\eeq
\begin{lem} \label{Levi-Civita}
The operator 
$\partial$ determines an isomorphism between the space of sections $\Sigma(E^* \otimes\so(E, g_o))$ and the space of sections  $\Sigma(\Lambda^2 E^* \otimes E)$. \end{lem}
\begin{pf} 
It suffices to check that, for any $x \in M$, the linear map
$\partial_x: ÊE^*_x\otimes\so(E_x, g_o) \longrightarrow \Lambda^2 E^*_x \otimes E_x$, defined by 
$
\partial H_x(s,s') = H_x(s, s') -  H_x(s', s)$, 
is an isomorphism. By dimension counting, it suffices to check that $\ker \partial_x = 0$. Following  a very classical argument,   this is proved noticing  that if $\partial H_x = 0$, then for any $s, s', s''$ one has 
$$g_o(H_x(s,s'), s'') = - g_o(H_x(s, s''), s') = - g_o(H_x(s'', s), s') = $$
$$ =  g_o(H_x(s'', s'), s) = 
g_o(H_x(s', s''), s) = - g_o(H_x(s', s), s'') = $$
$$ - g_o(H_x(s, s'), s'')$$ 
and hence that $H_x = 0$ by nondegeneracy of $g_o$.
\end{pf}
If $T^o$ is the torsion of $\n^o$, the torsion $T$ of any other derivation $(L, \n )$, represented by the pair $(L, \d \n)$,  is given by
 \beq \label{definitionoftorsion?}T = T^o + \partial \d \n\ .\eeq
From Lemma \ref{Levi-Civita},  $\d \n$ (and hence $\n$) can be completely recovered from $T$ and  $T^o$ and  the  following correspondence  is one-to-one:
\beq \label{deltaL} (L, T) \overset{\gd}\longrightarrow (L, \n) = \left(L, \n^o + \partial^{-1}(T- T^o)(L(\cdot), \cdot)\right)\eeq
\begin{lem} The correspondence  \eqref{deltaL}  is independent of  $\n^o$ and $T^o$ and gives a  one-to-one correspondence between the set of sections  $(L, T) \in \Sigma\left(T^*M \otimes_M E\right)\times\Sigma\left( \Lambda^2 E^*  \otimes E\right)$, with $L$ regular,  and the connections in $\Conn^{\psi_o}(\wt E)$.
\end{lem}
\noindent{\it Proof.}\ Let $\n^o$ and $\n^{'o}$ two metric covariant derivations of $(E, g_o)$, with torsions $T^o$ and $T^{'o}$, respectively,  and $\d \n^o \= \n^{'o} -  \n^o$. By (\ref{definitionoftorsion?}),   $T^{'o} = T^o + \partial \d \n^o$ and the conclusion follows from 
$$\n^{'o} + \partial^{-1}(T- T^{'o})(L(\cdot), \cdot) = $$
$$ = \n^{o} + \d \n^o + \partial^{-1}(T- T^{o})(L(\cdot), \cdot) - 
\partial^{-1}(\partial \d \n^O)(L(\cdot), \cdot) = $$
$$ =  \n^{o} + \ \partial^{-1}(T- T^{o})(L(\cdot), \cdot)\ .\eqno\qed$$
\begin{rem} By Lemma  \ref{Levi-Civita}, for any given $L$ there is a unique $\n^{oL}$ with vanishing torsion $T^o$. Inserting $\n^{oL}$ in  (\ref{deltaL}), the expression simplifies into 
 \beq \label{deltaL1} (L, T) \overset{\gd}\longrightarrow (L, \n) = \left(L, \n^{oL} + \partial^{-1}(T)(L(\cdot), \cdot)\right)\ .\eeq
\end{rem}
\medskip
\begin{coord}  In the notation used in \S \ref{classicalgravityfields} for the  expressions in  coordinates,  the torsion $T$ and the curvature $R$ of $(L,\n)$ are of the form 
$$T = T_{ij}^k e^i_o \otimes e^j_o \otimes e_k^o \ , \quad T^k_{ij} = e^\mu_i \Gam k \mu j - e^\mu_j \Gam k \mu i$$
$$R = R_{ijk}^m e^i_o \otimes e^j_o \otimes e^k_o \otimes e_m^o \ , \quad R^m_{ijk} = e^\mu_i e^\nu_j(\partial_\nu \Gam m \mu k - \partial_\mu \Gam m \nu k - \Gam m \mu \ell \Gam \ell \nu k + \Gam m \nu \ell \Gam \ell \mu k)\ .$$
Also the map $\partial^{-1}$ can be determined explicitly. It is equal to  
$$\partial^{-1}( T_{ij}^k e^i_o \otimes e^j_o \otimes e_k^o) = \frac{1}{2} \left(T_{ij}^k + T_{jk}^i - T_{ik}^j\right) e^i_o \otimes e^j_o \otimes e_k^o\ ,$$
which is the usual formula for the so-called ``contorsion''.
\end{coord}
\bigskip
\subsubsection{Infinitesimal pseudo-translations}\hfill\par
\medskip
In all the following, the section $\psi_o \in \Sigma(\wt E)$ is considered fixed and  $\Conn(\wt E)^{\psi_o reg}$ is identified with the  regular pairs $(L, \n)$. We introduce now the notion of pseudo-translations. In  the next Theorem \ref{main3} it is shown that  they  correspond to the infinitesimal transformations by vector fields on $M$.  \par
\medskip
\begin{definition} Let $X$ be a vector field on $M$. We call {\it infinitesimal pseudo-translation associated with $X$\/} the map 
$$\tau^{(X)}: \Conn(\wt E)^{\psi_o reg} \to \Sigma\left(T^*M \otimes_M E\right)\times\Sigma\left( E^* \otimes\so(E, g_o)\right)$$
\beq   \tau^{(X)}(L, \n) \= ( L\left(\n^LX  + T^L_{X\cdot}\right), \partial^{-1}(\delta_X T)(L(\cdot), \cdot))\eeq
where  $T^L$ is the torsion of  $\n^L$ (see \eqref{corr})  and $\delta_X T \in \Sigma\left( E^* \otimes\so(E, g_o)\right)$ is defined by 
$$(\delta_X T)_{s, s' }= T_{L(X) T(s, s')} +  T_{(\n_{L^{-1}(s)} L(X)) s'} - T_{(\n_{L^{-1}(s')} L(X)) s} +$$
 $$ - ( \n_{L^{-1}(s)} T) _{s'  L(X)} + ( \n_{L^{-1}(s')} T )_{s L(X)} + R_{ s s' L(X)} + R_{L(X) s s'} + R_{s' L(X) s}\ .$$
Here $T$ and $R$ are the torsion and curvature of $(L, \n)$ and  $\n_Y T$ is the vertical tensor defined by  $(\n_Y T)_{Z W} \= \n_Y(T_{Z W}) - T_{\n_Y Z  W} - T_{Z \n_Y W}$.
\end{definition}
Any infinitesimal pseudo-translation can be considered as  a ``vector field'' on $\Conn(\wt E)^{\psi_o reg}$ with associated flow $\cT^{(X)}_t: \Conn(\wt E)^{\psi_o reg} \to \Conn(\wt E)^{\psi_o reg}$ defined by
$$\left.\frac{d \cT^{(X)}_t (L, \n)}{dt}\right|_{t = t_o} = \tau^{(X)}(\cT^{(X)}_{t_o}(L, \n))\ . $$
We call  $\cT^{(X)}_t$   {\it flow of pseudo-translations generated by $X$\/}.\par
\medskip
The following theorem collects the main properties of pseudo-translations. In particular, it shows that  the action of the flow of a  pseudo-translation $\tau^{(X)}$  on the element $(L, \n) \in \Conn(\wt E)^{\psi_o reg}$ 
induces an action on the corresponding gravitational fields $(g, \n^L)$  on $M$ which coincide with the action of the flow of $X$ on $M$. 
\medskip
\begin{theo}\label{main3} \hfill\par
\begin{itemize}
\item[i)] For any $X \in \gX(M)$, the flow $\cT^{(X)}_t$  commutes with $\Gau_{\psi_o}(\wt E)$.
\item[ii)] Let  $\imath_{\psi_o}: \Conn(\wt E)^{\psi_o reg} \to \Grav_{p,q}(M)$ the correspondence \eqref{corr}, $(L, \n)\in\Conn(\wt E)^{\psi_o reg}$  and  $(g, \n^L) = \imath_{\psi_o}(L, \n)$.  For any  $X \in \gX(M)$
\beq\label{deforming}  \left.\frac{d }{dt}\imath_{\psi_o}(\cT^{(X)}_t (L, \n))\right|_{t = 0} = \left(\cL_X g, \left.\frac{d }{dt}\Phi^X_t{}_*( \n^L)\right|_{t = 0} \right)\ .\eeq
\item[iii)]ÊThe infinitesimal pseudo-translations have a natural structure of (infinite-dimensional) Lie algebra isomorphic to the Lie algebra of vector fields $\gX(M)$.
\end{itemize}
\end{theo}
\begin{pf}  To check (i), it is first necessary to observe that if $h \in \Gau_{\psi_o}(\wt E)$, then by definitions and \eqref{def_rot}, the torsion and curvature of $(h_* L, h_* \n)$ are equal to $h_* T = (h \circ T)(h^{-1} (\cdot), h^{-1}(\cdot))$ and $h_* R = (h \circ T)(h^{-1} (\cdot), h^{-1}(\cdot), h^{-1}(\cdot))$, respectively. Using this, a straightforward computation implies the claim.\par
 For (ii), let $(g^{(t)}, \n^{(t)}) = \imath_{\psi_o}(\cT^{(X)}_t (L, \n))$ and denote by $T^{(t)}$ the torsion of $\n^{(t)}$. Since $g^{(0)} = g$, $\n^{(0)} = \n^L$ and $T^L = T^{(0)}$,  we only need to  show that $\left.\frac{d g^{(t)}}{dt}\right|_{t = 0} = \cL_X g$ and  $\left.\frac{d T^{(t)}}{dt}\right|_{t = 0} = \cL_X T^{(0)}$. These identities are consequence of  the definition of pseudo-translations, standard properties of metric connections and first Bianchi identities. In fact, they can be obtained  observing that, for any $Y, Z \in \gX(M)$, 
 $$\cL_X g(Y, Z) = X(g(Y,Z)) - g(\cL_X Y, Z) - g(Y, \cL_X Z) =$$
 $$ = g(\n^L_Y X, Z) +  g(T^L_{X Y}, Z) + g(Y, \n^L_Z X) + g(Y, T^L_{X Z}) = $$
 $$ = g_o\left(L\left(\n^L_Y X + T^L_{X Y}\right),  L(Z)\right) + g_o\left(L(Y), L\left(\n^L_Z X + T^L_{X Z}\right)\right)$$
and 
  $$\cL_X T^L_{Y, Z} = X(T^L_{Y,Z})- T^L_{\cL_X Y Z} - T^L_{Y \cL_X Z} =   (\n^L_X T^L)_{Y,Z} +  T^L_{\n^L_Y X Z}  +  T^L_{T^L_{X Y}  Z} + $$
  $$ + T^L_{Y  \n^L_Z X} + T^L_{Y  T^L_{X Z}} =   - (\n^L_Y  T^L)_{Z X Z}  + (\n^L_Z T^L)_{Y X} + R_{X Y Z} + R_{Z X Y } + R_{Y Z X } +$$
$$   
  +  T^L_{\n^L_Y X Z}  +  T^L_{T^L_{X Y}  Z}\ . $$
Claim (iii) follows immediately from (ii) if we set $[\tau^{X}, \tau^{X'}] \= \tau^{[X, X']}$ the Lie bracket between two infinitesimal pseudo-translations. 
\end{pf}
By the above theorem,  we get the last result, mentioned  in the Introduction:  {\it an action on  pairs $(L, \n) \in \Conn(\wt E)^{\psi_o reg}$  corresponds to an action on gravitational fields $(g, \n^L)$, whose Euler-Lagrange equations  satisfy the infinitesimal  Equivalence Principle if and only if it  is invariant under pseudo-translations ``on-shell'', i.e. at the points given by solutions of the Euler-Lagrange equations\/}. \par
It is also clear that if $(L, \n)$ is so that 
 $T = 0$,  by definitions and  first Bianchi identities,  one has 
\beq   \tau^{(X)}(L, \n) =  (\n L(X) , \partial^{-1}(0))\eeq
and hence the flow $\cT^{(X)}_t$ maps  $(L, \n)$  into other metric affine connections with vanishing torsion, i.e. the condition $T = 0$ is preserved by pseudo-translations, as it should be by the correspondence   \eqref{deforming}.\par
\medskip
\subsubsection{A classical  example of gauge-invariant Lagrangian preserved by  pseudo-translations: the 
Palatini action}\hfill\par
Consider a space-time $M$ of dimension $n$,  a metric affine bundle $(\wt E, E, g_o)$ over $M$ with signature $(p,q)$,  and assume that there exists an affine vertical tensor field $\wt \o_o$, determined by a vertical volume form $\o_o$ on each fiber $E_x$, which is equal to $1$ on suitably ordered orthonormal frames $(e^o_j)$ of $E_x$.  A volume form $\o_o$ of this kind exists if and only if the bundle $O_{g_o}(E)$ of the orthonormal frames of  the fibers of $E$ admits an $SO_{p,q}(\bR)$-reduction.\par
Given  $(L, \n) \in \Conn(\wt E)^{\psi_o reg}$, we denote by 
 $\wt R$  the section in  $\Sigma(\Lambda^2 T^* M \otimes_M \Lambda^2 E)$  determined by the relation
$$e^i_o \otimes e^j_o(\wt R_{XY})  \= g_o(R_{L(X) L(Y)} \cdot e_i^o, e_j^o)$$
where $R$ is the curvature defined in Definition \ref{deftor} and $(e^i_o)$ is the coframe field dual to $(e_i^o)$.  Observe that, for any vector fields $Y_i \in TM$, $1 \leq i \leq n$ 
$$\o_o(\wt R_{Y_1 Y_2} \wedge L(Y_3) \wedge \dots \wedge L(Y_n))  = \sum_{i =1}^n \o_o (R_{Y_1 Y_2} \cdot e_i^o, 
e_i^o, L(Y_3), \dots, L(Y_n))\ .$$
We also denote by $\operatorname{Alt}: \otimes^n T^* M \to \Lambda^n T^*M$ the  usual alternating map. Now,  we consider the  {\it Palatini action\/}  on 
pairs  $(L, \n) \in \Conn(\wt E)^{\psi_o reg}$
$$\cS_{Pal}(L, \n) \= \int_M \operatorname{Alt}\left(\omega_o(\wt R \wedge L \wedge  \dots \wedge  L)\right) = $$
\beq  \label{palatini} 
 = \int_M \sum_{i =1}^n  \operatorname{Alt}\left(\o_o ((R \cdot e_i^o) \wedge  
e_i^o\wedge L \wedge \dots\wedge L\right) .  \eeq
 Using \eqref{def_rot},  one can check that the Lagrangian $\cL = \operatorname{Alt}\left(\omega_o(\wt R \wedge L \wedge  \dots \wedge  L)\right)$ is invariant under any rotation in $\cG_{\psi_o}$. If desired, one can also  extend 
 $\cL$ and obtain a fully gauge-invariant  action on $\Omega$ (see \S \ref{tggt}) by simply imposing that $\cL$ is  constant   along the  $\Gau^T(\wt E)$-orbits passing through the points of $\Conn(\wt E)^{\psi_o reg}$.  \par
The reader can also check that, expressed in terms of the pairs $(g^L, \n^L) \in \Grav_{p,q}(M)$,  the action $S_{Pal}$ becomes the usual Hilbert action $S_{Hilb} = \int_M \operatorname{Scal}(g) \o_g$, 
while the Euler-Lagrange equations determined by \eqref{palatini} coincide with those obtained  from 
the Palatini action with the Palatini method of variation (see e.g. \cite{Wi}), i.e. 
$$T^L = 0\ ,\qquad Ric^L = 0\ .$$ 
From Theorem \ref{main3}, it follows immediately that  \eqref{palatini} is also invariant under any pseudo-translations.\par
\begin{rem} In a future paper, we will consider gauge theories of  super-extensions of  Poincar\`e group  and the corresponding  analogues of pseudo-translations.  The  correspondence   with vector fields on  super-manifolds are expected to  relates the  invariance under pseudo-translations to a ``super''  version of the Equivalence Principle.
\end{rem}

\end{document}